\documentclass[centertags]{article}

\usepackage{amsmath}
\usepackage{amssymb}
\usepackage{latexsym}
\usepackage{amsxtra}
\usepackage{amscd}
\usepackage{theorem}
\usepackage{xypic}
\usepackage{epic}
\usepackage{eepic}

\theoremstyle{plain}

{\theorembodyfont{\slshape}

        \newtheorem{thm}{Theorem}[section]
        
        \newtheorem{lem}[thm]{Lemma}
        \newtheorem{prop}[thm]{Proposition}

}

{\theorembodyfont{\rmfamily}

        \newtheorem{defn}[thm]{Definition}
        
        \newtheorem{rem}[thm]{Remark}

        \newtheorem{exa}[thm]{Example}

}

\renewcommand{\em}{\sl}

\newcommand{\proof}{{\bf Proof:\ }}

\newcommand{\Endproof}{\hspace*{\fill} $\Box$ \vspace{1ex} \noindent }

\makeatletter
\renewcommand{\subsection}{\@startsection{subsection}{2}%
        {\z@}{-3.25ex plus -1ex minus-.2ex}{-1em}{\bf}}
\makeatother

\newcommand{\ZZ}{\mathbb{Z}}
\newcommand{\CC}{\mathbb{C}}

\newcommand{\QQ}{\mathbb{Q}}

\newcommand{\T}{\mathcal{T}}

\newcommand{\X}{\mathcal{X}}
\newcommand{\Y}{\mathcal{Y}}

\newcommand{\Kb}{\bar{K}}

\newcommand{\Xs}{\mathsf{X}}
\newcommand{\Ys}{\mathsf{Y}}

\newcommand{\Ds}{\mathsf{D}}
\newcommand{\Es}{\mathsf{E}}

\newcommand{\abs}[1]{\lvert#1\rvert}

\newcommand{\Aut}{\mathop{\rm Aut}\nolimits}

\newcommand{\Ker}{\mathop{\rm Ker}}

\newcommand{\val}{\mathop{\rm val}\nolimits}

\newcommand{\ord}{{\rm ord}}

\newcommand{\inj}{\hookrightarrow}
\newcommand{\To}{\longrightarrow}

\newcommand{\gen}[1]{\mathopen\langle#1\mathclose\rangle}

\newcommand{\Spf}{\mathop{\rm Spf}\nolimits}

\newcommand{\dcup}{\stackrel{\cdot}{\cup}}

\newcommand{\cyclic}[1]{\mathbb{Z}/#1\mathbb{Z}}

\title{Artin characters, Hurwitz trees\\ and the lifting
  problem} 
\author{Louis Brewis and Stefan Wewers}

\date{}

\begin{document}

\maketitle

\begin{abstract}
  We study finite groups of automorphisms of the $p$-adic open disk. In
  particular, we generalize results of Green, Matignon and Henrio from
  cyclic groups of order $p$ to arbitrary finite groups. As an application, we
  produce a counterexample to a question of Chinburg, Guralnick and
  Harbater, concerning the local lifting problem for generalized quaternion
  groups.
\end{abstract}

\section{Introduction}

\subsection{Automorphisms of the $p$-adic disk}

Let $K$ be a nonarchimedian local field of characteristic zero, with residue
field of characteristic $p$. Let $\Ys=(\Spf R[[z]])\otimes K$ be the rigid
open unit disk over $K$. In this note we are concerned with the following
question:
\begin{center}
  How can a finite group $G$ act on $\Ys$?
\end{center}
Let $G\subset\Aut_K(\Ys)$ be a finite group of automorphisms of $\Ys$;
alternatively, we may regard $G$ as a subgroup of $\Aut_R(R[[z]])$. It is not
hard to show that $G=P\rtimes C$ is the semi-direct product of a $p$-group $P$
and a cyclic group $C$ of order prime to $p$. Moreover, after a suitable
change of parameter the action of a generator $\sigma$ of $C$ is of the form
\[
      \sigma(z) = \zeta\cdot z,
\]
where $\zeta\in K^\times$ is a root of unity. In particular, the case where
$G$ is of order prime to $p$ is well understood.

In their fundamental paper \cite{GreenMatignon99}, Green and Matignon studied
the case of a cyclic group $G$ of order $p$ -- the first interesting
case. Their results establish a precise relation between
\begin{itemize}
\item the `geometry of fixed points', i.e.\ the relative position inside $\Ys$
  of the fixed points of a generator of $G$, and
\item
  the ramification of the Galois cover
  \[
       \Ys \To \Xs:=\Ys/G
  \]
  with respect to the Gauss norm of an arbitrary subdisk of $\Ys$.
\end{itemize}
These results were sharpened by Henrio (\cite{HenrioCR}, \cite{HenrioArxiv}),
in two ways. Firstly, he introduced the notion of {\em Hurwitz tree}, an
essentially combinatorial object which encodes both the geometry of fixed
points of an automorphism of order $p$ and the associated ramification
data. Secondly, he proved that every Hurwitz tree can be realized by
an actual automorphism of $\Ys$. 

\vspace{1ex} In the present paper we generalize the notion of Hurwitz tree to
the case of an arbitrary group of automorphisms of $\Ys$. Following
\cite{HenrioArxiv}, we define a Hurwitz tree as an oriented tree equipped with
a metric and some further additional data, satisfying certain conditions. Our
main results associates to a group of automorphisms $G$ of the disk a Hurwitz
tree $\T$. By construction, the metric tree underlying $\T$ describes the
geometry of fixed points of the $G$-action, and the additional data reflect
the ramification of the $G$-action with respect to certain valuations of
$R[[z]]$.  In the special case considered in \cite{HenrioArxiv}, these
additional data consist of certain numbers attached to each vertex and edge of
the tree. Our generalization is based on replacing these numbers by virtual
characters of $G$, which are constructed using Huber's ramification theory
\cite{Huber01}.

In this paper we do not consider the problem of realizing a given Hurwitz tree
by a $G$-action of the disk (which is solved in \cite{HenrioArxiv} for
$G\cong\ZZ/p$ and in \cite{Dp} for $G\cong\ZZ/p\rtimes\ZZ/m$, with
$(m,p)=1$). This explains why our Hurwitz trees are not equipped with {\em
  differential data}, as in \cite{HenrioArxiv} and \cite{Dp}. For a general
theory of Hurwitz trees with diffential data, see the first author's upcoming
thesis \cite{louisdiss}.

\subsection{New obstructions for the local lifting problem}

Let $k$ be an algebraically closed field of positive characteristic $p$. Let
$G$ be a finite group. A {\em local $G$-action} is a faithful $k$-linear
action $\phi:G\inj\Aut_k(k[[z]])$ of $G$ on a ring of power series over $k$ in
one variable. The {\em local lifting problem} for $\phi$ concerns the
following question:
\begin{quotation}
  Does $\phi$ lift to an action $\phi_R:G\inj R[[z]]$, where $R$ is a discrete
  valuation ring of characteristic $0$ with residue field $k$?
\end{quotation}
If it does, we say that {\em $\phi$ lifts to characteristic $0$}, or that
$\phi$ is {\em liftable}. 

In \cite{Bertin98}, Bertin has formulated a necessary condition for
liftability of $\phi$, which is commonly called the {\em Bertin
  obstruction}. If this necessary condition holds then we say that the Bertin
obstruction for $\phi$ {\em vanishes}. From our main result we directly obtain
a new necessary condition for liftability of local $G$-actions, which is a
refinement of the Bertin obstruction. 

We give examples of local $G$-actions (the so-called {\em simple
  quaternion actions}) with vanishing Bertin obstruction, but for which our
new condition does not hold -- and which therefore do not lift to
characteristic zero. Here the group $G$ is a generalized quaternion group
$Q_{2^{n+1}}$ of order $2^{n+1}$, with $n\geq 2$, and $k$ has characteristic
$2$.

Chinburg, Guralnick and Harbater (\cite{CGH1}, \cite{CGH2}) call a group $G$ a
{\em local Bertin group} if the Bertin obstruction of every local $G$-action
vanishes. They call $G$ a {\em local Oort group} if every local $G$-action
lifts to characteristic zero. They prove that the generalized quaternion
groups $Q_{2^{n+1}}$ are local Bertin groups for $n\geq 3$. However, our
result shows that these groups are not local Oort groups. This answers
Question 1.3 of \cite{CGH2} negatively.

\section{Hurwitz trees} \label{hurwitztree}

A Hurwitz tree $\T$ consists of an oriented metric tree $T$ and certain
additional data attached to each vertex and edge of $T$, satisfying certain
conditions. These additional data are related to a finite group $G$ and a prime
number $p$. Precise definitions will be given in the following Sections
\ref{htree1}, \ref{htree2} and \ref{htree3}. We postpone all motivation and
explanation of these definitions to Section \ref{disk}.

In Section \ref{density} we discuss the notion of {\em density}. Later on in
Section \ref{obstruction} this will be
our main tool for showing that certain Hurwitz trees and, therefore, certain
group actions on the disk, are impossible.

\subsection{} \label{htree1}

Let $G$ be a finite group. We denote by $R(G)$ the Grothendieck group of the
category of $\CC[G]$-modules of finite type (\cite{SerreRL}, \S 14.1).  We may
identify elements of $R(G)$ with their virtual characters $\chi:G\to\CC$. We
denote by $R^+(G)\subset R(G)$ the submonoid of true characters.

We write $1_G\in R^+(G)$ for the unit character, $r_G\in R^+(G)$ for the
regular character and $u_G=r_G-1_G\in R^+(G)$ for the augmentation character.

Given two characters $\chi_1,\chi_2\in R(G)$, their scalar product is defined
as 
\[
     \gen{\chi_1,\chi_2}_G := 
       \abs{G}^{-1}\cdot \sum_{\sigma\in G}\,
           \overline{\chi_1(\sigma)}\chi_2(\sigma).
\]
In the following, we will mostly identify an element of $R(G)$ with the
induced linear form on $R(G)$, i.e.\ we consider $\chi\in R(G)$ as a map
\[
     R(G)\to\ZZ, \qquad \psi\mapsto \chi(\psi):=\gen{\psi,\chi}_G.
\]

Given a group homomorphism $\phi:H\to G$, we obtain $\ZZ$-linear maps ({\em
  restriction} and {\em induction})
\[
       \phi^*:R(G)\to R(H), \qquad \phi_*:R(H)\to R(G).
\]
They are related by the {\em Frobenius reciprocity formula}:
\[
    \gen{\psi,\phi_*(\chi)}_G = \gen{\phi^*(\psi),\chi}_H.
\]
If $H$ is a subgroup of $G$ and $\phi$ the canonical injection, then we will
simply write $\psi|_H$ instead of $\phi^*(\psi)$ and $\chi^*$ instead of
$\phi_*\chi$. See e.g.\ \cite{SerreCL}, VI, \S 1.

A $\QQ$-valued virtual character is an element of the group
$R(G,\QQ):=R(G)\otimes_\ZZ\QQ$. We will consider an element $\chi\in R(G,\QQ)$
either as a class function $\chi:G\to\CC$ or as a $\ZZ$-linear map
$\chi:R(G)\to\QQ$. The submonoid $R^+(G,\QQ)\subset R(G,\QQ)$ consists, by
definition, of the elements $\chi\in R(G,\QQ)$
with $\chi(\psi)\geq 0$ for all $\psi\in R^+(G)$. An important example of such
a $\QQ$-valued character is the following:

\begin{defn} \label{deltamultdef} Let $p$ be a prime number and
  $G=\gen{\sigma}\cong\ZZ/p^m\ZZ$ be a finite cyclic group of order $p^m$,
  with $m\geq 0$. We define an element $\delta_G^{\rm mult}\in R(G,\QQ)$ via
  the following class function. For $a\not\equiv 0 \pmod{p^m}$ we set
 \[
     \delta_G^{\rm mult}(\sigma^a) := -\frac{p^{i+1}}{p-1},
 \]
 where $i:=\ord_p(a)<m$ is the exponent of $p$ in $a$; furthermore,
 \[
     \delta_G^{\rm mult}(1) :=  -\sum_{a=1}^{p^m-1}\,\delta_G^{\rm mult}(\sigma^a) 
                          = m\,p^m.
 \]
\end{defn}
 
Let $\chi\in R^+(G)$ be an irreducible character of $G$ of order $p^n$ (with
$0\leq n\leq m$). One easily checks that
\begin{equation} \label{deltamulteq}
     \delta_G^{\rm mult}(\chi) = \begin{cases}
       \;\;\frac{np-n+1}{p-1}, & n>0, \\
       \;\; 0,                 & \chi=1_G.
               \end{cases}
\end{equation}
It follows that $\delta_G^{\rm mult}\in R^+(G,\QQ)$. 

(The superscript $^{\rm mult}$ stands for {\em multiplicative} and was chosen
because $\delta^{\rm mult}$ describes the ramification of a torsor under the
multiplicative group scheme $\mu_{p^n}$. See Lemma \ref{mainthmlem}.)

\subsection{Metric trees} \label{htree2}

\begin{defn} \label{treedef}
  Let $T$ be a connected tree, with set of vertices $V$ and set of edges
  $E$ and with one distinguished vertex $v_0\in V$, called the {\em root}. We
  call $T$ a {\em rooted tree} if the root $v_0$ is connected to a unique edge
  $e_0\in E$ (which we call the {\em trunk} of $T$).
\end{defn}

A rooted tree $T$ carries a natural orientation, determined by source and
target maps $s,t:E\to V$, as follows. Given an edge $e\in E$, the source
$s(e)$ (resp.\ the target $t(e)$) is the vertex adjacent to $e$ contained in
same connected component of $T\backslash\{e\}$ as $v_0$ (resp.\ in the
connected component not containing $v_0$). If $v=s(e)$ and $v'=t(e)$ we call
$v'$ a {\em successor} of $v$; notation: $v\to v'$. There is a natural partial
ordering $\leq$ on $V$, where $v_1\leq v_2$ if and only if there is an
oriented path starting from $v_1$ and ending at $v_2$.

It is clear that the root $v_0$ is the unique minimal vertex with respect to
this ordering. A maximal vertex is called a {\em leaf}. We write $B\subset V$
for the set of all leaves. It follows from Definition \ref{treedef} that $B$
is nonempty and does not contain the root $v_0$.  For any vertex $v$ we define
\[
     B_v:=\{\,b\in B \,\mid\, v\leq b \,\}
\]
as the set of leaves which can be reached from $v$ along an oriented path.

\begin{defn}
  Let $T$ be a rooted tree. A {\em metric} on $T$ is given by a map
  $\epsilon:E\to\QQ_{\geq 0}$, $e\mapsto \epsilon_e$ such that $\epsilon_e=0$
  if and only if $t(e)$ is a leaf. 
  We call $\epsilon_e$ the {\em thickness} of the edge $e$.  The pair
  $(T,\epsilon)$ is called a {\em metric tree}. Sometimes we write $T$
  instead of $(T,\epsilon)$, if no confusion can arise.
\end{defn}

\subsection{Hurwitz trees} \label{htree3}

Let $G$ be a finite group and $p$ a prime number. 

\begin{defn} \label{def2.2} 
  A {\em Hurwitz tree} of type $(G,p)$ is a datum
  $\T=(T,[G_v],a_e,\delta_v)$, where
\begin{itemize}
\item
  $T=(T,\epsilon)$ is a metric tree (with root $v_0$, trunk $e_0$ and set of
  leaves $B$),
\item
  $[G_v]$ is the conjugacy class of a subgroup $G_v\subset G$, for every
  vertex $v$ of $T$, 
\item
  $a_e\in R^+(G)$ is a character of $G$, for every edge $e$ of $T$,
\item
  $\delta_v\in R^+(G,\QQ)$ is a $\QQ$-valued character of $G$, for all
  vertices $v$.
\end{itemize}
We call $G_v$ the {\em monodromy group} and $\delta_v$ the {\em depth} of the
vertex $v$. We call $a_e$ the {\em Artin character} of the edge $e\in E$.

The datum $\T$ is required to satisfy the following conditions:
\begin{itemize}
\item[(H1)]
  Let $v$ be a vertex. Then, up to conjugation in
  $G$, we have
  \[
       G_{v'}\subset G_v,
  \]
  for every successor $v'$ of $v$. Moreover, we have
  \[
          \sum_{v\to v'}\,[G_v:G_{v'}] > 1,
  \]
  except if $v=v_0$ is the root, in which case there exists exactly one
  successor $v'$ and we have $G_v=G_{v'}=G$.
\item[(H2)] 
  The group $G_b$ is nontrivial and
  cyclic, for every leaf $b\in B$.
\item[(H3)]
  For all $e\in E$ we have 
  \[
     a_e=\; \begin{cases}
            \;\;\sum_{t(e)=s(e')} a_{e'}, & t(e)\not\in B,\\
            \;\;u_{G_b}^*,         & b=t(e)\in B.
          \end{cases}
  \]
\item[(H4)]
  For all $e\in E$ we have
  \[
     \delta_{t(e)} = \delta_{s(e)}+\epsilon_e\cdot s_e,
  \]
  where $s_e:=a_e-u_{G_{t(e)}}^*\in R(G)$. 
\item[(H5)] For $b\in B$ we let $P_b\subset G_b$ denote the Sylow $p$-subgroup
  of $G_b$. Then
  \[
     \delta_b = (\delta_{P_b}^{\rm mult})^*.
  \]
  Here $\delta_{P_b}^{\rm mult}$ is given by Definition \ref{deltamultdef}
  (and depends on the prime $p$).
\end{itemize}
We set
\[
    \delta_\T:=\delta_{v_0}, \quad a_\T:=a_{e_0}, 
\]
which we call the {\em depth} and the {\em Artin character} of the Hurwitz
tree $\T$.
\end{defn}

\begin{rem} \label{rem2.2} 
  Let $\T=(T,[G_v],a_e,\delta_v)$ be a Hurwitz tree,
  as in Definition \ref{def2.2}.
\begin{enumerate}
\item
  Using Condition (H1), it is easy to show that there exists a metric tree
  $\tilde{T}$, together with an action of $G$ which fixes the root, such that 
  \begin{itemize}
  \item
    $T=\tilde{T}/G$, 
  \item
    for every vertex $v$ of $T$ there exists a vertex $\tilde{v}$ of
    $\tilde{T}$ above $v$ such that $G_v$ is the stabilizer of $\tilde{v}$. 
  \end{itemize} 
  In fact, a metric tree $\tilde{T}$ with $G$-action
  as above corresponds one-to-one to the datum $(T,G_v)$, satifying condition
  (H1). 
\item 
  Condition (H3) is equivalent to the following claim: for all edges $e$
  we have
  \[
          a_e= \sum_{b\in B_{t(e)}}\,u_{G_b}^*.
  \]
  This follows immediately from induction over the tree $T$. 
\item It follows from (ii) that the Artin characters $a_e$ are already
  determined by the tree $T$ and the conjugacy classes of (cyclic) subgroups
  $([G_b])_{b\in B}$. Moreover, using (H4) and (H5) we see that the depth
  $\delta_v$ is determined by the metrized tree $(T,\epsilon)$ and the
  conjugacy classes $([G_v])_{v\in V}$. 
\item By (iii) the Hurwitz tree $\T$ is uniquely determined by the datum
  $(T,[G_v])$ which corresponds, by (i), to a metric tree $\tilde{T}$ with a
  $G$-action. However, it is a highly nontrivial condition, for a metric tree
  $\tilde{T}$ with $G$-action, to come
  from a Hurwitz tree $\T$ of type $(G,p)$. 
\end{enumerate}
\end{rem}

\subsection{Densities} \label{density}

We fix a Hurwitz tree $\T=(T,[G_v],a_e,\delta_v)$ of type $(G,p)$, with set of
leaves $B$.

\begin{defn} \label{densitydef}
\begin{enumerate}
\item Let $b_1,b_2\in B$ be two distinct leaves. The {\em inverse distance}
  of $b_1$ and $b_2$ is the positive rational number
  $d(b_1,b_2)\in\QQ_{>0}$ defined as follows. Let
  $(v_0,v_1,\ldots,v_r)$ be the longest oriented path in $T$ starting
  from the root $v_0$ and ending in a vertex $v_r\not\in B$ with
  $v_r\leq b_1,b_2$. For $i=1,\ldots,r$ let $e_i$ be the edge with
  $s(e_i)=v_{i-1}$ and $t(e_i)=v_i$. Then we set
  \[
      d(b_1,b_2) := \sum_{i=1}^r \epsilon_{e_i}.
  \]
\item
  Let $A\subset B$ be a
  nonempty set of leaves and $b\in A$. The {\em density of $A$
  at $b$} is the rational number
  \[
     d(A,b) := \sum_{b'\in A\backslash\{b\}} d(b,b').
  \]
\end{enumerate}
\end{defn}
Note that $d(A,b)$ only depends on $A,b$ and the metrized tree $T$. 

\begin{lem} \label{densitylem}
  Let $A,b$ be as in Definition \ref{densitydef}. 
  \begin{enumerate}
  \item Let $(v_0,v_1,\ldots,v_r,b)$ be the unique oriented path
    from the root to $b$. For $i=1,\ldots,r$ let
    $e_i$ be the edge with $s(e_i)=v_{i-1}$ and $t(e_i)=v_i$. Then 
    \[
          d(A,b) = \sum_{i=1}^r \epsilon_{e_i}\cdot n(A,v_i),
    \]
    where
    \[
            n(A,v) := \abs{\{b'\in A \mid b'\neq b,\;v\leq b'\,\}}.
    \]
  \item
    Let $\chi\in R(G)^+$ be a character such that 
    \[
        \gen{\chi,u_{G_a}^*}_G = \begin{cases}
                \;\;m, & a\in A, \\
                \;\;0, & a\in B\backslash A,
             \end{cases}    
    \]
    where $m:=\gen{\chi,u_G}_G$. Then 
    \[
        m\cdot d(A,b) = \delta_b(\chi)-\delta_{v_0}(\chi).
    \]
  \end{enumerate}
\end{lem}

\proof The proof of (i) follows from a simple induction argument which
we leave to the reader. For the proof of (ii) we may assume that $G_b\subset
G_{v_r}\subset G_{v_{r-1}}\subset\ldots\subset G$, by Condition (H1) of
Definition \ref{def2.2}. We deduce the following sequence of inequalities
\[
     m=\gen{\chi,u_{G_b}^*}\leq \gen{\chi,u_{G_{v_r}}^*} \leq\ldots
       \leq \gen{\chi,u_G}=m,
\]
which, a posteriori, turn out to be equalities. Using Remark \ref{rem2.2} and
the hypothesis on $\chi$ we
therefore get 
\begin{equation} \label{lemdenseq1}
   s_{e_i}(\chi) = a_{e_i}(\chi) - m
     = \sum_{a\in B_{v_i}} \gen{\chi,u_{G_a}^*} -m 
                   = m\cdot n(A,v_i).
\end{equation}
Now we compute:
\[\begin{split}
    \delta_b(\chi)-\delta_{v_0}(\chi) 
      & =\quad \sum_{i=1}^r \delta_{v_i}(\chi)-\delta_{v_{i-1}}(\chi) \\
      & \stackrel{\rm (H3)}{=}\quad \sum_{i=1}^r \epsilon_{e_i}\cdot s_{e_i}(\chi) \\
      & \stackrel{\eqref{lemdenseq1}}{=}\quad  m\cdot \sum_i
      \epsilon_{e_i}\cdot n(A,v_i) \\
      & \stackrel{\rm (i)}{=}\quad m\cdot d(A,b).
\end{split}\]
\Endproof
 
\begin{exa} \label{CyclicExample} 
  Assume that $G = \cyclic{p^n}$ and let $b \in B$ such that $G_b = G$. Let
  $\chi$ be an irreducible character of $G$ with trivial kernel. Then by
  Lemma \ref{densitylem}, (H5) and \eqref{deltamulteq} we have
   \[
      d(B,b) = \frac{np - n + 1}{p-1}-\delta_\T(\chi).
   \]
   If $\delta_\T=0$ (which is the interesting case for us) we thus get
   a simple formula for the density $d(B,b)$ which puts a strong restriction
   on the metric of the tree $T$.
\end{exa}




\section{Group actions on the disk} \label{disk} 

\subsection{}

We fix the following notation. Let $K$ be a field of characteristic zero which
is complete with respect to a discrete (additive) valuation
$\val:K^\times\to\QQ$. Let $R\subset K$ denote the valuation ring and $\pi$ a
prime element of $R$. We assume that the residue field $k=R/(\pi)$ is
algebraically closed of characteristic $p>0$. We also assume that $\val(p)=1$.

We fix an open rigid-analytic disk $\Ys$ over $K$ and a subgroup
$G\subset\Aut_K(\Ys)$. We assume that there exists at least one fixed point,
i.e.\ a point in $\Ys$ with a nontrivial stabilizer. The goal of this section
is to attach to $(\Ys,G)$ a Hurwitz tree $\T=(T,[G_v],a_e,\delta_v)$ of type
$(G,p)$.

This construction is based on Huber's theory of Artin and Swan characters for
rigid-analytic curves (\cite{Huber01}). But since we only consider a very
special case (a disk), we can do everything in an elementary and
self-contained way, and we do not have to actually use any of the results of
\cite{Huber01}.

\subsection{The depth character}

At the beginning we shall work with a slightly more general situation than
announced above. The ring $A$ will either denote the ring of formal power
series $R[[z]]$ or the ring $R\{z\}$ of convergent powers series in $z$. It
gives rise to a formal scheme $\Y:=\Spf A$ and a rigid-analytic space
$\Ys:=\Y\otimes K$. In the first case, $\Ys$ is an {\em open disk}, i.e.
\[
    \Ys = \{\,z \,\mid\, \val(z) >0 \,\}.
\]
In the second case it is a {\em closed disk}, and we have a bijection
\[
    \Ys = \{\, z \,\mid\, \val(z) \geq 0 \,\}.
\]
We let $\val_\Ys:A\backslash\{0\}\to\QQ$ denote the Gauss valuation, i.e.\
\[
     \val_\Ys\big(\,\sum a_iz^i\,\big) = \min_i\,\val(a_i).
\]
We set $\bar{A}:=A/(\pi)$ and let $\bar{f}\in\bar{A}$ denote the image of $f\in
A$. We have $\bar{A}=k[[z]]$ or $\bar{A}=k[z]$. 

Suppose we are given a finite subgroup $G\subset\Aut_K(\Ys)$ of automorphisms
of $\Ys$. The action of $G$ extends uniquely to the formal model $\Y$ and
hence induces an action of $G$ on the ring $A$. 

Our first goal is to define an invariant $\delta_\Ys^G\in R^+(G,\QQ)$, called
the {\em depth character}. It measures the ramification of $G$ with respect to
$\val_\Ys$, i.e.\ the amount to which the induced map $G\to\Aut_k(\bar{A})$
fails to be injective. 

Let $I\lhd G$ be the inertia group with respect to $\val_\Ys$, i.e.\ the
normal subgroup consisting of elements $\sigma\in G$ with
$\val_\Ys(\sigma(z)-z)>0$.

\begin{defn} \label{depthdef}
  The {\em depth character} associated to $(\Ys,G)$ is the $\QQ$-valued
  character $\delta_\Ys^G\in R(G,\QQ)$ associated to the following class
  function: 
  \[
       \delta_\Ys^G(\sigma) := -\abs{G}\cdot\val_\Ys(\sigma(z)-z)
  \]
  for $\sigma\in G\backslash\{1\}$ and
  \[
       \delta_\Ys^G(1) := -\sum_{\sigma\neq 1}\,\delta_\Ys^G(\sigma).
  \]
\end{defn}
By definition we have $\delta_\Ys^G=0$ if and only if $I=\{1\}$. 

A priori, it is not clear why $\delta_\Ys^G$ is a class function, as it seems
to depend on the choice of the parameter $z$. However, the following lemma
shows that it actually does not.

\begin{lem} \label{depthlem1}
  Fix an element $\sigma\in I\backslash\{1\}$. Then:
  \begin{enumerate}
  \item
    For all $f\in A$ we have
    \[
       \val_\Ys(\sigma(f)-f) \geq \val_\Ys(\sigma(z)-z).
    \]
  \item
    We have equality in (i) if and only if ${\rm d}\bar{f}\neq 0$.
  \end{enumerate}
  Therefore, $\delta_\Ys^G(\sigma)$ is independent of
  the choice of the para\-meter $z$ and only depends on the conjugacy class of
  $\sigma$ in $G$.
\end{lem}

\proof
Since $\sigma\in I\backslash\{1\}$, there exists an element $a\in R$ such that
\[
     \val(a) = \min\{ \val_\Ys(\sigma(f)-f) \mid f\in A \}>0.
\]
For $g\in A$ with $\val_\Ys(g)\geq \val(a)$ we define
\[
     [g] := \overline{(g/a)}\in\bar{A}.
\]
We claim that the map
\[
     \partial:\bar{A}\to\bar{A}, \qquad \bar{f}\mapsto [\sigma(f)-f].
\]
is a derivation. Indeed, it is clearly $k$-linear,
and for $f,g\in A$ we get:
\[\begin{split}
   \partial(\overline{fg}) & = [\sigma(f)(\sigma(g)-g)+g(\sigma(f)-f)] \\
     &= \bar{f}\,\partial(\bar{g}) +\bar{g}\,\partial(\bar{f}) 
         + [(\sigma(f)-f)(\sigma(g)-g)] \\
     &= \bar{f}\,\partial(\bar{g}) +\bar{g}\,\partial(\bar{f}).
\end{split}\]
This proves the claim. But $\partial\neq 0$ by definition. It follows
that $\partial(\bar{f})=0$ if and only if ${\rm d}\bar{f}=0$. This proves the
lemma. 
\Endproof

\begin{lem} \label{depthlem2}
    We have $\delta_\Ys^G\in R^+(G,\QQ)$.
\end{lem}

\proof One easily checks (using Lemma \ref{depthlem1}) that for $h\in\QQ_{\geq
  0}$ the set
\[
    G_h:= \{\sigma\in G \mid \val_\Ys(\sigma(z)-z)\geq h \}
\]
is a normal subgroup of $G$. Let
\[
    0<h_1<\ldots< h_m
\]
be the breaks, i.e.\ the positive values of the function
\[
    G\backslash\{1\} \to \QQ, \qquad \sigma\mapsto\val_\Ys(\sigma(z)-z).
\]
One verifies, by a direct computation, that
\[
    \delta_\Ys^G = \sum_{i=1}^m \lambda_i\cdot u_{G_{h_i}}^*,
\]
where
\[
    \lambda_i:= \abs{G_{h_i}}(h_i-h_{i-1}) >0.
\]
and $h_0:=0$. The lemma follows from the positivity of $\lambda_i$.
\Endproof

\subsection{The Artin character}

We continue with the notation introduced above. But from now on we assume that
$A=R[[z]]$, i.e.\ that $\Ys$ is an open disk. Our goal is to define an Artin
character $a_\Ys^G\in R^+(G)$ which describes the action of $G$ on the
boundary of $\Ys$.

We define
\[
    \#_\Ys f := \ord_z\big(\overline{f/p^{\val_\Ys(f)}}\big).
\]
Here $\ord_z:k[[z]]\to\ZZ\cup\{\infty\}$ is the usual order function and
$p^{\val_\Ys(f)}\in R$ is an arbitrary element with valuation $\val_\Ys(f)$.
The Weierstrass preparation theorem shows that $\#_\Ys f$ is the number of
zeroes of $f$ on $\Ys$, counted with multiplicity.

\begin{defn} \label{artindef}
  The {\em Artin character} of $(\Ys,G)$ is the element of $R^+(G)$ associated
  to the class function defined by 
  \[
     a_\Ys^G(\sigma) := -\#_\Ys (\sigma(z)-z), \qquad 
         \text{for $\sigma\neq 1$}
  \]
  and 
  \[
      a_\Ys^G(1):=-\sum_{\sigma\neq 1}\,a_\Ys^G(\sigma).
  \]
\end{defn}

To see that $a_\Ys^G$ is indeed a character we relate it to the permutation
representation arising from the set of fixed points. For $\sigma\in
G\backslash\{1\}$ let $\Delta_\sigma\subset\Ys(\Kb)$ denote the set of
(geometric) fixed points of $\sigma$. Set
\[
    \Delta:=\cup_{\sigma\neq 1} \,\Delta_\sigma.
\]
This is a finite $G$-set. Let $B:=\Delta/G$ denote the orbit space. Choose,
for each $b\in B$, an element $y\in\Delta$ belonging to $b$ and let
$G_b\subset G$ denote the stabilizer of $y$. 

\begin{prop} \label{artinprop}
  We have 
  \[
      a_\Ys^G = \sum_{b\in B} \, u_{G_b}^*.
  \]
  In particular, $a_\Ys^G$ is an element of $R^+(G)$.
\end{prop}

\proof Fix an element $\sigma\in G\backslash\{1\}$. Then $\Delta_\sigma$ is
the set of zeroes of the function $f_\sigma:=\sigma(z)-z$. An easy local
calculation, coupled with the assumption that $\sigma$ has finite order and
that ${\rm char}(K)=0$, shows that all zeroes of $f_\sigma$ are simple (cf.\
\cite{GreenMatignon99}, \S II.1). Therefore, by Definition \ref{artindef} and
the Weierstrass preparation theorem we have
\[
      a_\Ys^G(\sigma) = -\#_\Ys f_\sigma = -\abs{\Delta_\sigma}.
\]
The proposition follows immediately.
\Endproof

The next proposition is the key result behind the construction of the
Hurwitz tree associated to $(\Ys,G)$. 

\begin{prop} \label{keyprop} Let $\Ds\subset\Ys$ be a closed disk which
  contains the set $\Delta$ and is fixed by the action of $G$. Let
  $\Es\subset\Ds$ denote the residue class of a $K$-rational point $y$ in
  $\Ds$. Let $H\subset G$ denote the stabilizer of $\Es$. Then
  \begin{equation} \label{keypropeq}
     \delta_\Ds^G = \big(\delta_\Es^H)^* = 
       \delta_\Ys^G + \abs{G}\cdot\epsilon\cdot s_\Ys^G,
  \end{equation}
  where $s_\Ys^G:=a_\Ys^G-u_G$ and where $\epsilon\in\QQ_{>0}$ is the
  thickness of the annulus $\Ys\backslash\Ds$. 
\end{prop}

\proof After a change of parameter we may assume that the point $y$ is given
by the equation $z=0$. Then 
\[
     \Ds = \{\, z \mid \val(z)\geq \epsilon\,\}, \quad
     \Es = \{\, z  \,\mid\, \val(z)>\epsilon \,\}.
\]
After replacing $K$ by some finite extension, we may further assume that there
exists an element $a\in R$ with $\val(a)=\epsilon$. We obtain formal models
$\Ds=(\Spf R\{w\})\otimes_RK$ and $\Es=(\Spf R[[w]])\otimes K$, where
$w:=a^{-1}z$. By definition, we have $\val_\Ds=\val_\Es|_{R\{w\}}$ and
therefore
\[
    \delta_\Ds^G(\sigma) = 
     \begin{cases}
       \;\;[G:H]\cdot\delta_\Es^H(\sigma), & \sigma\in H\backslash\{1\} \\
        \qquad 0, & \sigma\in G\backslash H.
     \end{cases}
\]
Now the first equation in \eqref{keypropeq} is obvious.

Fix an element $\sigma\in G\backslash\{1\}$ and set $f_\sigma:=\sigma(z)-z\in
R[[z]]\subset R\{w\}$. By the assumption on $\Ds$, the function $f_\sigma$ has
no zero on the annulus $\Ys\backslash\Ds$. It follows that
\begin{equation} \label{valeq}
  \val_\Ds(f_\sigma) = \val_\Ys(f_\sigma) + \epsilon\cdot \#_\Ys f_\sigma,
\end{equation}
see e.g.\ \cite{Amice}, \S 4.5.  We compute:
\[\begin{split}
  \delta_\Ds^G(\sigma) &=\; -\abs{G}\cdot\val_\Ds(\sigma(w)-w) 
     = -\abs{G}\cdot \big(\val_\Ds(f_\sigma)-\epsilon\big) \\
     &\stackrel{\eqref{valeq}}{=}\; -\abs{G}\cdot \val_\Ys(f_\sigma) 
          - \abs{G}\cdot\epsilon\cdot(\#_\Ys(f_\sigma)-1) \\
     &=\; \delta_\Ys^G(\sigma) + \abs{G}\cdot\epsilon\cdot s_\Ys^G.
\end{split}
\]
This proves the second equation in \eqref{keypropeq}. 
\Endproof

\subsection{Definition of the Hurwitz tree}

We can now state and prove our main theorem.

\begin{thm} \label{mainthm}
  Let $\Ys=(\Spf R[[z]])\otimes K$ be an open rigid disk over $K$ and
  $G\subset\Aut_K(\Ys)$ be a finite group of automorphisms. Suppose that the
  set of fixed points $\Delta\subset \Ys$ is nonempty. Then there exists a
  Hurwitz tree $\T$ of type $(G,p)$ with 
  \begin{equation} \label{mainthmeq}
      \delta_\T=\delta_\Ys^G, \qquad a_\T = a_\Ys^G.
  \end{equation}
\end{thm}

\proof Our proof is by induction over the number of elements of $\Delta$. 

\noindent
We first assume that $\abs{\Delta}=1$. In this case the theorem is essentially
equivalent to the following lemma. 

\begin{lem} \label{mainthmlem}
  Let $y\in\Delta$ be the unique fixed point. Then 
  \begin{enumerate}
  \item 
    the group $G$ is cyclic,
  \item
    $a_\Ys^G=u_G$, and
  \item
    $\delta_\Ys^G=(\delta_P^{\rm mult})^*$, where $P$ is the Sylow
    $p$-subgroup of $G$. 
  \end{enumerate}
\end{lem}

\proof It is clear that every element of $G$ fixes the point $y$. So (ii)
follows directly from Proposition \ref{artinprop}. 

After a change of parameter we may assume that $y$ is the point $z=0$. Then
for an element $\sigma\in G$ we have
\begin{equation} \label{mainthmlemeq1}
     \sigma(z) =\chi(\sigma)\, z\,(1+a_1z+a_2z^2+\ldots),
\end{equation}
where $\chi:G\inj K^\times$ is an injective character (\cite{GreenMatignon99},
\S II.1). This proves (i). Let us fix an element $\sigma\in G$ of order
$np^m$, with $(n,p)=1$ and $m\geq 0$. By \eqref{mainthmlemeq1} we have
\[
   f_\sigma:=\sigma(z)-z=(\chi(\sigma)-1)\,z +\chi(\sigma)\,a_1z^2+\ldots
\]
Since $z=0$ is the only zero of $f_\sigma$, we have $\#f_\sigma=1$ and therefore
\[
   \val_\Ys(f_\sigma) = \val(\chi(\sigma)-1) =
     \begin{cases}
        \;\;0, & m=0, \\
        \;\;\frac{1}{(p-1)p^{m-1}}, & m>0.
     \end{cases} 
\]
Now (iii) follows from Definition \ref{deltamultdef} and a direct
computation. 
\Endproof

So in the case $\abs{\Delta}=1$ we define the Hurwitz tree
$\T=(T,[G_v],\delta_v,a_e)$ as follows.
\begin{itemize}
\item
  The tree $T$ has two vertices $v_0,v_1$ and one edge $e_0$ with $s(e_0)=v_0$
  and $t(e_0)=v_1$. The metric $\epsilon$ is trivial, i.e.\ we set
  $\epsilon_{e_0}:=0$. 
\item
  We define 
  \[
     \delta_{v_0}=\delta_{v_1}:=\delta_\Ys^G.
  \]
\item
  We define $G_{v_0}=G_{v_1}:=G$ and $a_{e_0}:=u_G$.
\end{itemize}
The validity of the axioms (H2) and (H5) follows from Lemma \ref{mainthmlem};
all the other axioms and \eqref{mainthmeq} hold by definition. This finishes
the proof of the theorem in the case $\abs{\Delta}=1$.

We may now assume that $\abs{\Delta}\geq 2$. Then there exists a smallest
closed disk $\Ds\subset\Ys$ which contains $\Delta$. Clearly, $\Ds$ is fixed
by the $G$-action. There also exists a finite family $(\Es_j)_{j\in
  J}$ of residue classes $\Es_j\subset \Ds$ with
\begin{equation} \label{mainthmeq0}
   \Delta_j:=\Es_j\cap\Delta\neq\emptyset \quad\text{and}\quad
    \Delta\subset \cup_j\Es_j.
\end{equation}
For $j\in J$ we let $G_j\subset G$ denote the stabilizer of $\Es_j$. By
induction, there exists a Hurwitz tree $\T_j$ for the group $G_j$ with
\begin{equation} \label{mainthmeq1}
     \delta_{\T_j} = \delta_{\Es_j}^{G_j}, \qquad
     a_{\T_j} = a_{\Es_j}^{G_j}.
\end{equation}

The Hurwitz tree $\T=(T,[G_v],a_e,\delta_v)$ associated to $(\Ys,G)$ is defined
as follows.
\begin{itemize}
\item Let $T_j$ denote the metric tree underlying the Hurwitz tree
  $\T_j$. Choose a system of representatives $J'\subset J$ of $J/G$.  The
  metric tree $T$ underlying $\T$ is obtained by patching together the metric
  trees $T_j$, $j\in J'$, at their roots, i.e.\ we identify the set of roots
  of the trees $T_j$, $j\in J'$ with one vertex $v_1$ of $T$. We complete $T$
  by adding another vertex $v_0$ (the root of $T$) and an edge $e_0$ with
  $s(e_0)=v_0$, $t(e_0)=v_1$. The value of the metric $\epsilon$ on the edge
  $e_0$ is defined as the thickness of the annulus $\Ys\backslash\Ds$,
  multiplied with $\abs{G}$.  (In fact, $\epsilon_{e_0}$ is the thickness of
  the quotient annulus $(\Ys\backslash\Ds)/G$.)
\item
  If $v$ is a vertex of $T$ other than $v_0$ and
  $v_1$, it corresponds to a vertex $v'$ of one of the $T_j$ which is not the
  root. We define $G_v:=G_{v'}$ and $\delta_v:=\delta_{v'}^*$.
\item
  Let $e$ be an edge of $T$ which corresponds to an edge $e'$ of $T_j$. We
  define $a_e:=a_{e'}^*$.
\item
  We set $G_{v_0}=G_{v_1}:=G$, $\delta_{v_0}:=\delta_\Ys^G$,
  $\delta_{v_1}:=\delta_\Ds^G$ and $a_{e_0}:=a_\Ys^G$.   
\end{itemize} 
It remains to show that $\T$ satisfies the axioms (H1)-(H5). Since these
axioms hold for the Hurwitz trees $\T_j$, many of them hold for $\T$ by
construction. For instance, this is clear for (H1) and (H2).

It follows from \eqref{mainthmeq0}, \eqref{mainthmeq1} and Proposition
\ref{artinprop} that
\begin{equation} \label{mainthmeq3}
     a_{e_0}=a_\Ys^G = \sum_{j\in J/G} \, (a_{\Es_j}^{G_j})^*=\sum_{s(e)=v_1}a_e.
\end{equation}
Therefore, (H3) holds for the edge $e_0$. For the other edges it holds by
construction. 

To check the axioms (H4) and (H5) we remark that 
\begin{equation} \label{mainthmeq2}
     \delta_{v_1}=\delta_\Ds^G = \delta_{\T_j}^*,
\end{equation}
for all $j\in J$, by the first equality in \eqref{keypropeq}. This means that
our definition of $\delta_{v_1}$ is consistent with the fact that the vertex
$v_1$ corresponds to the roots of the Hurwitz trees $\T_j$, $j\in J'$. It
follows that (H5) holds automatically and that we have to check (H4) only for
the edge $e_0$. But for the edge $e_0$ the statement of (H4) follows directly
from Proposition \ref{keyprop}.  This concludes the proof of Theorem
\ref{mainthm}.
\Endproof

\begin{rem}
  An alternative way to construct the metric tree $T$ is the following (cf.\
  \cite{HenrioArxiv} and \cite{Dp}). Let $\Y$ be the minimal semistable model
  of the disk $\Ys$ which separates the points of $\Delta$. Then the
  $G$-action on $\Ys$ extends to $\Y$, and the quotient $\X:=\Y/G$ is a
  semistable model of the disk $\Xs=\Ys/G$ which separates the points of
  $B:=\Delta/G$. Now there is a standard way to associate to the pair $(\X,B)$ a
  metric tree $T$ with set of leaves $B$ (see e.g.\ \cite{Dp}, \S
  3.2). Essentially, $T$ is a modification of the graph of components of the
  special fiber of $\X$.

  The construction of $T$ in the proof of Theorem \ref{mainthm} avoids the use
  of semistable models and may therefore be considered as more
  elementary. However, semistable models become inevitable if one wants to
  construct $G$-actions on the disk with given Hurwitz tree.
\end{rem}




\section{Applications to the lifting problem} \label{obstruction}

\subsection{A new obstruction}

Let $k$ be an algebraically closed field of characteristic $p>0$ and $G$ be a
finite group. A {\em local $G$-action} is a faithful and $k$-linear action
$\phi:G\inj\Aut_k(k[[z]])$ on a ring of formal power series in one variable
over $k$. 

The {\em local lifting problem} asks: can $\phi$ be lifted to an action
$\phi_R:G\inj\Aut_R(R[[z]])$, where $R$ is some discrete valuation ring of
characteristic zero with residue field $k$. If it does then we say that $\phi$
{\em lifts to characteristic zero}.

From our main result we can deduce a new necessary condition for liftability
of local $G$-actions. Before we state it, we recall the definition of the {\em
  Artin character}.  

\begin{defn}
  Let $\phi$ be a local $G$-action. The {\em Artin character} of $\phi$ is the
  element $a_\phi\in R^+(G)$ defined by
  \[
      a_\phi(\sigma) := -\ord_z(\sigma(z)-z)
  \]
  for $\sigma\neq 1$ and 
  \[
     a_\phi(1) := -\sum_{\sigma\neq 1}\,a_\phi(\sigma).
  \]
\end{defn}
See \cite{SerreCL}, VI, \S 2.

\begin{thm} \label{obstrthm}
  Let $\phi:G\inj\Aut_k(k((t)))$ be a local $G$-action. If $\phi$
  lifts to characteristic $0$ then there exists a Hurwitz tree
  $\T$ of type $(G,p)$ such that
  \[
       a_\T = a_\phi \qquad\text{and}\qquad \delta_\T=0.
  \]
\end{thm}

\proof A lift of $\phi$ gives rise to a $G$-action on the disk $\Ys=(\Spf
R[[z]])\otimes K$. Since $\phi$ is injective by assumption, we have
$\delta_\Ys^G=0$ (Definition \ref{depthdef}) and $a_\Ys^G=a_\phi$ (Definition
\ref{artindef}). Therefore, Theorem \ref{obstrthm} is a direct consequence of
Theorem \ref{mainthm}. 
\Endproof 

By the theorem, the existence of a Hurwitz tree $\T$ with given
Artin character $a_{\T}=a_\phi$ and trivial depth $\delta_{\T}=0$ is a
necessary condition for $\phi$ to lift. If one can show that such a Hurwitz
tree does not exist, one has found an obstruction against liftability of
$\phi$. 

As a special case of this criterion, we obtain the well-known {\em Bertin
  obstruction}, see \cite{Bertin98}. Namely, if $\T=(T,[G_v],a_e,\delta_v)$ is a
Hurwitz tree with $a_\T=a_\phi$, then Remark \ref{rem2.2} shows that
\begin{equation} \label{bertineq1}
    a_\phi = \sum_{b\in B} u_{G_b}^*.
\end{equation}
This equality is easily seen to imply the following statement: there exists a
finite $G$-set $\Delta$, with cyclic stabilizers, such that
\begin{equation} \label{bertineq}
         a_\phi = m\cdot r_G - \chi_\Delta. 
\end{equation}
Here $\chi_\Delta\in R^+(G)$ is the character of the permutation
representation realized by $\Delta$ and $m:=\abs{\Delta/G}$. However, there
exist local $G$-actions $\phi$ whose Artin character can {\em not} be written
in this form (see e.g.\ \cite{Bertin98} and \cite{CGH2}). It follows from
Theorem \ref{obstrthm} that such a $\phi$ does not lift to characteristic
zero.

\begin{rem}
  The examples presented in Section \ref{simple} show that our new obstruction
  is strictly stronger than the Bertin obstruction. However, it should be
  pointed out that the converse of Theorem \ref{obstrthm} does not hold. For
  $G=\ZZ/p\times\ZZ/p$, Pagot has shown in \cite{Pagot} that certain local
  $G$-actions $\phi$ do not lift to characteristic zero. For such a $\phi$ it
  is straightforward to write down a Hurwitz tree $\T$ with $a_\T=a_\phi$ and
  $\delta_\T=0$.
\end{rem}

\begin{rem}
  It will be shown in \cite{louisdiss} show that our new obstruction vanishes
  for all cyclic groups $G$, in accordance with Oort's conjecture (see
  \cite{CGH2}). The interesting thing about the proof is that one has to use
  certain nontrivial inequalities satisfied by the Artin character of a local
  action of a cyclic group of order $p^n$ (see \cite{HLSchmidt37}). We believe
  that this observation provides further substantial evidence in favour of
  Oort's conjecture.
\end{rem}

\subsection{Simple quaternion actions} \label{simple}

We fix an integer $n\geq 2$ and let $G=Q_{2^{n+1}}$ denote the
  generalized quaternion group of order $2^{n+1}$, with presentation
\[
     G=\gen{\, \sigma,\tau \,\mid\, 
         \tau^{2^n}=1,\, \tau^{2^{n-1}}=\sigma^2,\, 
           \sigma\tau\sigma^{-1}=\tau^{-1} \,}.
\]
Our base field $k$ is assumed to be of characteristic $2$.

Chinburg, Guralnick and Harbater (\cite{CGH2}) have proved that $G$ is a {\em
  local Bertin group} for $n\geq 3$, which means that the Bertin obstruction
of every local $G$-action over $k$ vanishes. The goal of this section is to
construct certain $G$-actions which do not lift to characteristic zero. This
result gives a negative answer to Question 1.3 of \cite{CGH2}. 

We first introduce some more notation. Set
\begin{equation} \label{cyclicgrouplist}
     H_0:=\gen{\tau}, \quad H_1:=\gen{\sigma},\quad H_2:=\gen{\sigma\tau};
\end{equation}
these are cyclic subgroups of $G$ of order $2^n$, $4$ and $4$,
respectively. For $i=0,1,2$ there exists a unique character $\chi_i:G\to\{\pm
1\}$ of order $2$ such that $H_i\subset\Ker(\chi_i)$. Clearly,
$\chi_0,\chi_1,\chi_2$ define pairwise distinct irreducible characters of the
quotient group
\[
     \bar{G} := G/\gen{\tau^2} \cong \ZZ/2\times\ZZ/2.
\]
The following lemma is an easy exercise:

\begin{lem} \label{cycliclem}
  If $H\subset G$ is a cyclic subgroup whose image in $\bar{G}$ is nontrivial,
  then $H$ is conjugate to one of the groups in \eqref{cyclicgrouplist}.
\end{lem}

\begin{defn} \label{simpledef}
  A local $G$-action $\phi$ is called {\em simple} if 
  \[
     a_\phi(\chi_0) =2, \qquad a_\phi(\chi_1)=a_\phi(\chi_2)\geq 2.
  \]
\end{defn}

\begin{prop}
  There exists a simple $G$-action over $k$, for every $n\geq 2$.
\end{prop}

\proof
Choose an embedding of abelian groups $\bar{G}\inj(k,+)$. We obtain a local
$\bar{G}$-action $\bar{\phi}:\bar{G}\inj\Aut_k(k[[t]])$ by sending
$\mu\in\bar{G}$ to the automorphism
\[
    t\mapsto \frac{t}{1+\mu t} = t-\mu t^2+\mu^2 t^3-\ldots
\]
One checks that 
\[
   a_{\bar{\phi}}(\chi_i)=2, \qquad \text{for $i=0,1,2$.}
\]
By \cite{CGH1}, Lemma 2.10, we can extend $\bar{\phi}$ to a local $G$-action
$\phi:G\inj\Aut_k(k[[z]])$, such that $k[[t]]=k[[z]]^{\gen{\tau^2}}$. It follows
from \cite{SerreCL}, Proposition IV.3, that
\[
     a_\phi(\chi_i)=a_{\bar{\phi}}(\chi_i), \qquad i=0,1,2.
\]
We conclude that $\phi$ is simple. 
\Endproof

\begin{thm} \label{simplethm} Let $\phi$ be a simple $G$-action over $k$. Then
  $\phi$ does not lift to characteristic zero.
\end{thm}

\proof Suppose that $\phi$ lifts to characteristic zero. By Theorem
\ref{obstrthm}, there exists a Hurwitz tree $\T=(T,[G_v],a_e,\delta_v)$ of type
$(G,2)$, with Artin character $a_\T=a_\phi$ and vanishing depth
$\delta_\T=0$. We will show that such a Hurwitz tree cannot exist. Our main
tool is the notion of density introduced in Section \ref{density}.

Let $B$ denote the set of ends of the tree $T$. For $i=0,1,2$ we set
\[
    B_i:=\{b\in B \mid [G_b]=[H_i] \}, \qquad B':= B_0\dcup B_1\dcup B_2
\]
and 
\[
     B^i:=B'\backslash B_i. 
\]
It follows from Lemma \ref{cycliclem} and the definition of the groups $H_i$
that 
\begin{equation} \label{simpleeq2}
    \gen{\chi_i,u_{G_b}^*} = \; 
       \begin{cases}
         \;\; 1, & b\in B^i, \\
         \;\; 0, & b\in B\backslash B^i,
       \end{cases}
\end{equation}
for all $b\in B$. So Lemma \ref{densitylem} (ii), Condition (H5) of Definition
\ref{def2.2} and \eqref{deltamulteq} show that
\begin{equation} \label{simpleeq3}
  d(B^i,b) = \delta_b(\chi_i)= 2,
\end{equation}
for all $b\in B^i$. 

From \eqref{bertineq1} and \eqref{simpleeq2} we conclude that
\begin{equation} \label{simpleeq4}
  a_\phi(\chi_i) = \abs{B^i} = \sum_{j\neq i} \abs{B_j}.
\end{equation}
Now the assumption that $\phi$ is simple (Definition \ref{simpledef}) implies
that
\[
    \abs{B_1}=\abs{B_2}=1, \qquad \abs{B_0}\geq 1.
\]

Let $b_0$ denote the unique element of $B_2$. Since the $B_i$ are disjoint, we
have $B^0\cap B^1=B_2=\{b_0\}$ and $B^0\cup B^1=B'$. Using Definition
\ref{densitydef} (ii) and \eqref{simpleeq3} we therefore get
\begin{equation} \label{simpleeq5}
  d(B',b_0) = d(B^0,b_0) + d(B^1,b_0) = 2+ 2 = 4.
\end{equation}    

Let $\chi:H_0\inj\CC^\times$ be an injective irreducible character. The
induced character $\psi:=\chi^*\in R^+(G)$ has the following property. For any
nontrivial cyclic subgroup $C\subset G$, the restriction $\psi|_C$ is the sum
of two {\em nontrivial} irreducible characters of $C$,
$\psi|_C=\psi_1+\psi_2$. Applying this to $C=G_b$, we obtain
\[
      \gen{\psi,u_{G_b}^*}_G=\gen{\psi_1,u_{G_b}}+\gen{\psi_2,u_{G_b}} = 2,
\]
for all $b\in B$. We may therefore apply Lemma \ref{densitylem} (ii) and
conclude that 
\begin{equation} \label{simpleeq6}
  d(B,b_0) = \delta_{b_0}(\psi)/2. 
\end{equation}
Moreover, the restriction of $\psi$ to $G_{b_0}=H_2\cong\ZZ/4$ is the sum of
two irreducible characters $\psi_1,\psi_2$ of order $4$. From (H5) and
\eqref{deltamulteq} we get
\begin{equation} \label{simpleeq7}
  \delta_{b_0}(\psi)= \delta_{H_2}^{\rm mult}(\psi_1)+
         \delta_{H_2}^{\rm mult}(\psi_2) = 3+3=6.
\end{equation}
We now obtain a contradiction by comparing \eqref{simpleeq5},
\eqref{simpleeq6} and \eqref{simpleeq7}:
\[
   4=d(B',b_0)\leq d(B,b_0)= 3.
\]
We conclude that there does not exist a Hurwitz tree $\T$ of type $(G,2)$ with
$a_\T=a_\phi$ and $\delta_\T=0$. Theorem \ref{simplethm} follows.
\Endproof

\vspace{4ex}
{\small
\begin{minipage}[t]{5cm}
Institut f\"ur Reine Mathematik\\
Universit\"at Ulm\\
Helmholtzstr.\ 18\\
89069 Ulm\\
louis.brewis@uni-ulm.de\\
\end{minipage}
\hfill
\begin{minipage}[t]{5cm}
\begin{flushright}
IAZD\\
Leibniz-Universit\"at Hannover\\    
Welfengarten 1\\
30167 Hannover\\
wewers@math.uni-hannover.de
\end{flushright}
\end{minipage}
}

\end{document}